\documentclass[10pt,a4paper,twoside]{article}

\usepackage{amsfonts, amssymb, amsmath, amsthm}
\usepackage{mathrsfs} % \mathscr
\usepackage{latexsym}%utilisation des symboles LaTeX pour avoir un beau LaTeX
\usepackage{enumerate}
\usepackage{multicol}
\usepackage{verbatim}
\usepackage{dsfont}

\usepackage{url}

\newtheorem{thm}{Theorem}[section]
\newtheorem{lem}[thm]{Lemma}
\newtheorem{cor}[thm]{Corollary}

\newcommand{\Z}[1]{\mathbb{Z}/#1\mathbb{Z}}

\DeclareMathOperator{\ch}{ch}
\DeclareMathOperator{\sh}{sh}

\def\Acal{\mathcal{A}}
\def\Bcal{\mathcal{B}}
\def\Ccal{\mathcal{C}}

\def\Xcal{\mathcal{X}}

\def\Ocal{\mathcal{O}}

\def\1{\mathds{1}}

\def\mode{\mathbin{\,\textrm{mod}^*}}
\def\mod{\mathbin{\,\textrm{mod}\,}}

\title{Rudin Inequality, Chang Theorem, primes and squares%
}
\author{Olivier Ramar\'e}

\begin{document}
% \address[O. Ramar\'e]{CNRS/ Institut de Math\'ematiques de Marseille, Aix 
% Marseille Universit\'e, U.M.R. 7373, Site Sud, Campus de Luminy, Case 907, 
% 13288 
% Marseille Cedex 9, France.}
% \email{olivier.ramare@univ-amu.fr}

% %\date{\sl January, the 7th of 2004}
% \subjclass[2010]{Primary: , Secondary: }

% \keywords{Class Field}

\maketitle

\hfill \parbox{7.28cm}{\sl Dedicated to George Andrews and Bruce
  Berndt for their 85${}^{\text{th}}$ birthdays} 

\begin{abstract}
  %\texttt{File \jobname.tex}
  
We prove that the set of large values of
the trigonometric polynomial over a subset of density of the primes has some additive structure, similarly to what happens for
subsets of densities in $\Z{N}$ but in a weaker form.
To do so, we prove large sieve inequalities for \emph{dissociate sets}
$\mathcal{X}$ of circle 
points and functions $f$ whose support~$S$ is finite and
respectively in an interval, in
the set of primes or in the set of squares. Set
$T(f,x)=\sum_{n}f(n)\exp(2i\pi nx)$. These inequalities are of
the shape $\sum_{x\in\mathcal{X}}|T(f,x)|^2\ll
|S|\|f\|_2^2\log(8R/|S|)$ where $R$ is respectively $N$, $N/\log
N$ and $\sqrt{N}$. The implied constants depend on
the spacement between sumsets of~$\mathcal{X}$.
%We deduce them from related Rudin type inequalities.
\end{abstract}

%{\small \tableofcontents}

%%%%%%%%%%%%%%%%%%
\section{Introduction and results}
%%%%%%%%%%%%%%%%%%
One of the main outcome of the present paper is that the set of large values of
the trigonometric polynomial over a subset of density of the primes
has some additive structure, similarly to what happens for 
subsets of densities in $\Z{N}$. In the case of squares, which is
extreme with respect to the method employed, we reach only an
improved large sieve inequality.

Such a property has been proved for dense subsets of integers by
M.~Chang~\cite{Chang*02} as a consequence of an inequality attributed
to W.~Rudin, explaining our title; the setting generally employed is the
one of finite abelian groups,
see for instance the survey paper \cite{Green*04} by B.~Green or the
book \cite{Tao-Vu*06} by T. Tao \& V.H. Vu. We consider here
subsequences of the primes from $[1,N]$ and want to rely on sieve
techniques. These
rely crucially on the fact that the primes are invertible
(if we omit the initial ones) in $\Z{q}$ for any $q\le \sqrt{N}$; this
fact is not naturally expressible inside~$\Z{\Ocal(N)}$, so
we first establish some results on integer sequences
lying in~$[1,N]$. These are not new in spirit, though some constants
may be improved with respect to the published ones. The author does
not know of any ancestors to the case of
primes and of squares that is examined later.

A main definition on this query comes from harmonic analysis.
In accordance with the book
\cite[Definition 2.5]{Lopez-Ross*75} by J.~L\'opez \& K.~Ross, a finite set
$\mathcal{X}\subset \mathbb{R}/\mathbb{Z}$ is called \emph{dissociate} when the
equation $\sum_{x\in\mathcal{X}}\epsilon(x) x=0$ has no non-trivial
solution with $(\epsilon(x))\in\{0,\pm1\}^{\mathcal{X}}$; a natural
example is the set $\{\frac{1}{U},\frac{2}{U},\cdots,\frac{2^H}{U}\}$
for some positive integers $U$ and $H$. We present
more methological and historical points below. 
An inequality of Rudin involving such sets of points and that is usually
traced to~\cite{Rudin*60} plays an important role in our proof.
The paper
\cite{Shkredov*11} by I.~Shkredov has a very informative first part
with many references and may be used as an introduction to the
subject. We somehow displace the center of gravity of the questions at
hand towards an improved large sieve inequality over such sets. The
usual large sieve inequality says that (see \cite{Montgomery*78}
by H.L. Montgomery), for any finite subset 
$\mathcal{X}\subset\mathbb{R}/\mathbb{Z}$, we have
\begin{equation}
  \label{LSI}
  \sum_{x\in\mathcal{X}}\biggl|\sum_{n\le N}f(n)e(nx)\biggr|^2
  \le
  (N+\delta^{-1})\, \|f\|_2^2
\end{equation}
where $\|f\|_2^2=\sum_b |f(n)|^2$,  $e(y)=\exp(2i\pi y)$  and
\begin{equation}
  \label{eq:7}
  \delta=\min\bigr\{\|x-x'\|_{\mathbb{R}/\mathbb{Z}}:
  x,x'\in\mathcal{X},
  x\neq x'\bigl\}.
\end{equation}
Notation $\|u\|_{\mathbb{R}/\mathbb{Z}}$ stands for the distance on
the unit circle, i.e.
\begin{equation}
  \label{eq:6}
  \|u\|_{\mathbb{R}/\mathbb{Z}}=\min_{k\in\mathbb{Z}}|u-k|.
\end{equation}
In the case of applications, we have $\delta\gg 1/N$.
By assuming a quantitative separation condition on the points of
$\mathcal{X}$, namely that $\delta_\star\gg 1/N$, where $\delta_\star$
is defined in~\eqref{defdeltastar}, we prove the improved large sieve
like inequality given in Theorem~\ref{LpQRI*}. It is stronger in that,
in essence, we replace the $N+\delta^{-1}$ above by the length of the
support of~$f$.

When working over finite abelian groups, say over $\Z{U}$, it is
enough to use the condition that all the sumsets of $\mathcal{X}$ are
distincts. In the case of an interval, and since we do not restrict
our attention to circle points of the shape $u/U$,  the
quantitative information provided to us by $\delta_\star$ is needed.

For primes, we require a stronger and
more arithmetical quantity, namely $\delta_*(z,z_0)$ defined in
\eqref{defdeltastarp}. For intervals, the linear forms in
the points of $\mathcal{X}$ with
coefficients $\{0,\pm1\}$ had to be away from~0. We now need these
forms to be away from any rational $a/q$ with a not-too-large denominator.
Of course asking for $\delta_*(z,z_0)>0$ is
stronger than asking that $\delta_\star>0$, but applications ask in
fact for $\delta_*(z,z_0)\gg 1/N$ which is this time much stronger
than the condition $\delta_\star\gg 1/N$.
We show in Corollary~\ref{cuspsmodU} below how to treat this
requirement geometrically.
The parameter $z_0$ is included to allow
some small prime factors in~$U$. Having this parameter
$\delta_*(z,z_0)$ at hand, we prove two improved large sieve
inequalities that link dissociate sets on one side and subset of
density either of primes or of squares on the other side. Notice that
this means a two-fold saving: first we need to save the density of the
primes (or of the squares) with respect to the integers as in
\cite{Green-Tao*04} by B.~Green \& T.~Tao and more efficiently in
\cite{Ramare*22}, and then the 
relative density of the support of $f$ with respect to the primes (or
to the squares).

This is achieved by using an \emph{enveloping sieve}, as developed in
\cite{Ramare*91}, in~\cite{Ramare*95} and in \cite{Ramare-Ruzsa*01}
(see also \cite{Ramare*12-a} for a reminder). We further employ a
Fourier analytic device (Lemma~\ref{BeurlingF}) to reduce the arithmetical input to a
minimum (in the proof of Theorem~\ref{LpQRI*p} and~\ref{LpQRI*Q}, we
evaluate only a density and do not have to handle any error
term). The enveloping sieve for the squares is much less used, so we spend
some time putting it in place.
 
%%%%%%%%%%%%%%%%%%%%%%
\subsection*{The case of intervals}
%%%%%%%%%%%%%%%%%%%%%%
We proceed by producing a general enough proof than will prove at the
same time several Rudin's like inequalities, infer a large sieve
inequality from a dual version of it and deduce a Chang's theorem
from there. We first prove an interval version.
%%%%%%%%%%%%%%
\begin{thm}
  \label{LpQRI*}
  Let $\Xcal\subset\mathbb{R}/\mathbb{Z}$ be a finite set. Let
  \begin{equation}
    \label{defdeltastar}
    \delta_\star=\min
    \biggl\{\Bigl\|\sum_{x\in\Acal}x-\sum_{x\in\Bcal}x\Bigr\|_{\mathbb{R}/\mathbb{Z}}:
    \Acal\neq\Bcal\subset\Xcal\biggr\}.
  \end{equation}
  Assume that $\delta_\star>0$. When $f$ has support
  inside $S\subset\{1,\ldots,N\}$, we have
  \begin{equation*}
    \sum_{x\in\Xcal}
    \biggl|\sum_{n\le N}f(n)e(nx)\biggr|^2
    \le
    9
    \, |S|\, \|f\|_2^2\, \log\frac{8(N+\delta_\star^{-1})}{|S|}
    .
  \end{equation*}
  More generally, for any real number $\ell\ge1$, we have
  \begin{equation*}
    \biggl(\sum_{x\in\Xcal}
    \biggl|\sum_{n\le N}f(n)e(nx)\biggr|^{\ell+1}\biggr)^2
    \le
    9\,
    |S|\,\|f\|_2^2\,\sum_{x\in\Xcal}
    \biggl|\sum_{n\le N}f(n)e(nx)\biggr|^{2\ell}
    \log\frac{8(N+\delta_\star^{-1})}{|S|}
    .
  \end{equation*}
\end{thm}
%%%%%%%%%%%%%%
\noindent
Let us compare with the large lieve inequality~\eqref{LSI}: we
have here $|S|\log\frac{8N}{|S|}$ instead of $N$, on assuming that
$\delta_\star\gg1/N$. The dependence in~$\Xcal$ is however
much worse.
The parameter~$\delta_\star$ is absent from the theory developed on
finite abelian groups. See Corollary~\ref{cuspsmodU} for a result
without this parameter.
The corresponding Rudin's inequality is given in Lemma~\ref{LpQRI}.

I.~Shkredov establishes in \cite[Theorem 1.3]{Shkredov*11} the
analogous of the inequality with~$\ell$, save that he requires the
condition~$\ell\ge2$ that we can wave in our proof.
The $\ell$-part for $\ell>2$ may also be proved by using the case $\ell=2$: we
apply the Cauchy inequality to the scalar product
$\langle u|v\rangle=\sum_{x}|T(f,x)|^2u(x)\overline{v(x)}$ and the
variables $u=\1$ and $v=|T(f,x)|^{\ell-1}$, where
$T(f,x)=\sum_{n}f(n)\exp(2i\pi nx)$. Though we do not specify it, the
extension to $\ell\ge1$ is valid for
Theorem~\ref{LpQRI*p} and~\ref{LpQRI*Q} and for
Corollary~\ref{cuspsmodU}.

%%%%%%%%%%%%%%%%%%%%%%
\subsection*{The case of primes}
%%%%%%%%%%%%%%%%%%%%%%

Here our first main result, namely the counterpart of
Theorem~\ref{LpQRI*} for primes. 
%%%%%%%%%%%%
\begin{thm}
  \label{LpQRI*p}
  Let $\Xcal\subset\mathbb{R}/\mathbb{Z}$ be a finite set. Let
  $N\ge z_0\ge 2$ be two parameters and define $P(z_0)=\prod_{p<z_0}p$. We set
  \begin{equation}
    \label{defdeltastarp}
    \delta_*(z,z_0)
    =\min\biggl\{
    \Bigl\|\sum_{x\in\Acal}x-\sum_{x\in\Bcal}x-\frac{a}{q}\Bigr\|_{\mathbb{R}/\mathbb{Z}}:
    \Acal\neq\Bcal\subset\Xcal,
    a\in\mathbb{Z}, q\le z, (q,P(z_0))=1
    \biggr\}.
  \end{equation}
  Let $\kappa\in(0,1/2]$.
  Assume that $\delta_*(N^\kappa,z_0)>0$. For any function $f$ with support
  inside a subset~$S$ of the primes of $\{N^\kappa,\ldots,N\}$, we have
  \begin{equation*}
    \sum_{x\in \Xcal}\biggl|
    \sum_{p\le N}f(p) e(xp)\biggr|^2
    \le c(\kappa)|S|\,\|f\|_2^2\,
    \log\frac{8N}{|S|\log N}
    .
  \end{equation*}
  where
  \begin{equation*}
    c(\kappa) = 9\frac{\log\Bigl(64\frac{N+\delta_*^{-1}(N^\kappa,z_0)}{\kappa
        |S|\log N}\log
      z_0\Bigr)}{\log\frac{8N}{|S|\log N}}.
  \end{equation*}
\end{thm}
%%%%%%%%%%%%
\noindent
The corresponding inequality in \cite{Ramare*22} has essentially
$N(\log N)^{-1}\log |\Xcal|$ when the above with $z_0=2$ has
$N(\log N)^{-1}(\log K)/K$, with $K=N/(|S|\log N)$. On assuming that
$c$ is indeed bounded above (i.e. that
$\delta_*(\sqrt{N},z_0)\gg 1/N$), the saving is thus two-fold: we save
(almost all) the relative density $K$ of the subset~$|S|$ and
the~$\log|\Xcal|$.

This result is more easily read on rational circle points with a
denominator with a large prime divisor, where we
recover the group structure on one component.  Here is a sample
result.
%%%%%%%%%%%%
\begin{cor}
  \label{cuspsmodU}
  Let $N\ge2$ and let $U\ge N$ be an integer with a prime divisor
  $U_1\in[N^{1/4}, N^{3/4}]$. Let
  $\mathcal{U}\subset\Z{U}$ all whose sumsets modulo~$U_1$ are distinct.
  For any function $f$ with support
  inside a subset~$S$ of the primes of $\{N^{1/4},\ldots,N\}$, we have
  \begin{equation*}
    \sum_{u\in \mathcal{U}}\biggl|
    \sum_{p\le N}f(p) e\bigg(\frac{up}{U}\biggr)\biggr|^2
    \le 30\,|S|\,
    \|f\|_2^2\, e^{U/N}    \,\log\frac{8N}{|S|\log N}.
  \end{equation*}
\end{cor}
%%%%%%%%%%%%
\noindent
From the facts that the number of subsets of $\mathcal{U}$ is $2^{|\mathcal{U}|}$ and that
each such subset corresponds to only one sum modulo~$U_1$, we infer that
$|\mathcal{U}|\le (\log U_1)/\log 2$. This bound is optimal as the
example $\{1,2,\cdots,2^H\}$ shows. We may also replace $e^{U/N}$
by $e^{U/(N\log 8K)}$; we chose simplicity.

Here is now a results in the line of M.~H.~Chang in \cite[Lemma 3.1]{Chang*02},
but for dense prime subsets.
%%%%%%%%%%%
\begin{thm}[Chang's theorem for primes]
  \label{StrutCusps}
  Let $N\ge2$ and let $U_1$ and $U_2$ be two distinct primes such that
  $U_1,U_2\in[N^{1/4}, N^{3/4}]$. Set $U=U_1U_2$.
  Let $S$ be a subset of the primes of $\{N^{1/4},\ldots,N\}$ and let
  $A\ge1$ be given. Set 
  $K=N/(|S|\log N)$. For $i\in\{1,2\}$, there exists $\mathcal{D}_i\subset\Z{U_i}$ of
  cardinality at most $30A^2\log(8K)e^{U/N}$ such that
  \begin{equation*}
    \biggl\{u\in\Z{U}: \biggl|\sum_{p\in S}e(pu/U)\biggr|\ge
      |S|/A\biggr\}
    \subset
    \prod_{i=1}^2\biggl\{\sum_{a\in\mathcal{D}_i}\epsilon(a)a:
    (\epsilon(a))\in\{0,\pm1\}^{\mathcal{D}_i}\biggr\}.
  \end{equation*}
\end{thm}
%%%%%%%%%%%
\noindent
It is a consequence of the material exposed in~\cite{Ramare*22} that,
in case $U\ll N$,
the set on the left hand side, say $\mathscr{C}$, has cardinality at most
$\Ocal(A^2 K\log 2A)$. The above theorem
therefore exhibits some non-trivial additive structure of the
set~$\mathscr{C}$ when~$A$ is small with respect to~$K$. The
set $\mathcal{D}$ may indeed be smaller in size than $\mathscr{C}$, though 
the set of the
right-hand side above may have cardinality $9^{\Ocal(A^4\log^2
  (8K))}$.

Let us end this part with a counter-example that explains why we do
not restrict our attention to a dissociate set
$\mathcal{D}_1\times\mathcal{D}_2$. The projection of the dissociate
set 
$\mathcal{D}=\{(1,1),(1,4),(3,2),(3,3)\}\subset\Z{5}\times\Z{6}$
on $\Z{5}$ is the dissociate set $\{1,3\}$, but there are
not dissociate subset $\mathcal{A}$ of $\Z{6}$ such that
$\mathcal{D}\subset\{1,3\}\times \mathcal{A}$.  
%%%%%%%%%%%%%%%%%%%%%%
\subsection*{The case of squares}
%%%%%%%%%%%%%%%%%%%%%%
Theorem~\ref{LpQRI*p} can immediately be generalized to any sieve
context of finite dimension, as in~\cite{Green-Tao*04} by B.~Green \& T.~Tao.
As it turns out, the Selberg sieve has teeth even for sieving squares,
and still provides there an upper bound that is only a constant times
bigger than expected. This is not the case anymore if we sieve prime squares. All that
has been developped and used in \cite[Theorem 5.4]{Ramare*06}. 
It leads readily to the next theorem.
%%%%%%%%%%%%
\begin{thm}
  \label{LpQRI*Q}
  With hypothesis and notation as in Theorem~\ref{LpQRI*p} with $z_0=2$.
  For any function $f$ with support
  inside a subset~$S$ of the squares of $\{1,\ldots,N\}$, we have
  \begin{equation*}
    \sum_{x\in \Xcal}\biggl|
    \sum_{n\le \sqrt{N}}f(n) e(xn^2)\biggr|^2
    \le c|S|
    \|f\|_2^2\,\log\frac{8\sqrt{N}}{|S|}
    .
  \end{equation*}
  where the quantity $c$ is given by
  \begin{equation*}
    c =
    9\frac{\log\Bigl(64\frac{N+\delta_*^{-1}(\sqrt{N},2)}{|S|\sqrt{N}}\Bigr)}
    {\log\frac{8\sqrt{N}}{|S|}}.
  \end{equation*}
\end{thm}
%%%%%%%%%%%%
\noindent
The squares in the context of $\Lambda(p)$-sets have been intensively examined.
In this extreme case, we need to sieve truely up to $\sqrt{N}$, so we
cannot afford a parameter~$\kappa$ as in the case of primes, and that
ruins any corollary like Corollary~\ref{cuspsmodU}.

%%%%%%%%%%%%%%%%%%%%%%
\subsection*{Some methological and historical notes}
%%%%%%%%%%%%%%%%%%%%%%
%%%%%%%%%%%%%%%%%%%%%%
%\subsection*{Notes on $\Lambda(p)$-sets}
%%%%%%%%%%%%%%%%%%%%%%
Following \cite[Section 4.5]{Tao-Vu*06} of the book by T.~Tao \& V.~Vu,
we define the $\Lambda(p)$-constant of a subset $S\subset\Z{U}$ to be
the smallest $K$ such that
\begin{equation}
  \label{defLambdapsets}
  \biggl(\frac1U\sum_{u\mod U}\biggl|\sum_{s\in S}c(s)e(su/U)\biggr|^p\biggr)^{1/p}
  \le K
  \biggl(\frac1U\sum_{u\mod U}\biggl|\sum_{s\in S}c(s)e(nu/U)\biggr|^2\biggr)^{1/2}.
\end{equation}
This inequality should hold for any choice $(c(s))$; $S$ is a
subset of the dual group of the initial goup ($X$ in the notation of
the reference book~\cite{Lopez-Ross*75} by J.~L\'opez \& K.~Ross, the
initial group being denoted by~$G$).  Such an inequality is close to
what we call a generic Rudin's inequality in Lemma~\ref{LpQRI} below
and is generally traced back to~\cite[Theorem 3.1]{Rudin*60} by
W.~Rudin. In harmonic analysis over groups, see for instance
Definition 5.2 in~\cite{Lopez-Ross*75}, a subset $S$ is called a
$\Lambda(p)$-set if any such constant exists, but in the finite set
context, only an upper bound for the attached constant is meaningful.

When we replace $\Z{U}$ by the interval~$[1,N]$, it is clear that we
stay close to the definition~\eqref{defLambdapsets}, but the
counterparts with primes or squares are less readable in this
context. The other variable, which is from $\mathcal{X}$ in our
notation, comes naturally from the character group of $\mathbb{Z}$,
and may be also seen, when containing only rational points, as a
subset of the character group of some $\Z{U}$. This one corresponds
truely to the $s$-variable in \eqref{defLambdapsets}.

The notion of dissociate set appears only in the
1969-paper~\cite{Hewitt-Zuckerman*69} by E.~Hewitt \& H.~Zuckerman,
about ten years after the Rudin paper. However \cite[Theorem
2.4]{Rudin*60} contains points that shows that a dissociate set is a
\emph{Sidon} set.  Yet again, the notion of Sidon sets is somewhat
difficult in the context of finite structure. One definition used in
number theory is that no point of a Sidon set equals to a sum of two
points from this set. This is not the notion used harmonic analysis,
but there are translations of the harmonic analytic definition in
arithmetical context, see for instance \cite[Theorem 2.4]{Rudin*60} or
the book \cite[Theorem 5.7.5]{Rudin*62} again by W. Rudin. Yet again,
this definition looses its impact in finite structures if one does not
follow the size of the implied finite parts.

The intermediate notion of \emph{$\Lambda(p)$-sets} seems more
appropriate in the context of finite sets, save for the caveat that
web search engines have difficulties with such a name! We refer to the
survey paper \cite{Bourgain*01} by J.~Bourgain. In this paper,
Eq.~(3.7) defines as an \emph{independent} set what appears to be a
dissociate set, bar the possibility of characters of order~2.  To add
to this joyful confusion, it seems modern terminology prefers
\emph{dissociated} to \emph{dissociate}.

As we have seen, the notion of dissociate sets which is instrumental in
our work appears in earlier work only on the character group side,
and the same holds for the notion of Sidon sets. There are deep
questions therein concerning squares, dissociate sets and Sidon sets. 
Our Theorem~\ref{LpQRI*Q} which mixes this kind of conditions on both
variables is thus rather new. 

%%%%%%%%%%%%%%%%%%%%%% 
\subsection*{Extensions and acknowledgements}
%%%%%%%%%%%%%%%%%%%%%%
In \cite{Bourgain*08f}, J.~Bourgain proves an extension to
Theorem~\ref{LpQRI*} to sumsets of a dissociated set. The reader may
read \cite{Shkredov*09} by I.~Shkredov for a proper introduction to
the subject. The possibility of such an extension to the case of
primes or of squares are open questions.

We decided to use a rather low-level style of writing so as to be
accessible to both communities of number theorists and harmonic
analysists. We hope this will be helpful to the readers.

A large part of this work was completed when the author was enjoying
the hospitality of the Indian Statistical Institute in Kolkata. It has
also has been partially supported by the Indo-French IRL Relax
and by the joint FWF-ANR project Arithrand: FWF: I 4945-N and ANR-20-CE91-0006.
The author also thanks warmly H.~Queff\'elec and I.~Shkredov for their fruitful comments.

%%%%%%%%%%%%%%%%%% 
\section{An abstract machinery}
%%%%%%%%%%%%%%%%%%
\label{Abstract}
We want to prove things on the integers and then mimick the proof to
reach corresponding results on the primes. On looking closely, one
sees that one input differs and the rest follows.

In this section, we consider a finite sequence $\mathcal{N}$ of
integers, a finite sequence~$\mathcal{X}$ of points of
$\mathbb{R}/\mathbb{Z}$. We assume that there exists a constant $H\ge |\mathcal{N}|$
such that, for every complex sequence $(c(x))_{x\in\mathcal{X}}$, we
have
\begin{equation}
  \tag{Hyp.}
  \log \sum_{n\in\mathcal{N}}\biggl|
  \exp\biggl(
  \sum_{x\in \Xcal}c(x) e(xn)\biggr)\biggr|
  \le
  \tfrac12 \sum_{x\in\Xcal}|c(x)|^2
  +
    \log H
    .
\end{equation}
%%%%%%%%%%%%%%%%%%
\paragraph{Remark:}
%%%%%%%%%%%%%%%%%%
On choosing $c(x)=0$ when $x\neq0$, for some chosen $x_0$ from
$\mathcal{X}$, the above inequality reduces to
$  -(\tfrac12-c(x_0))^2+\tfrac14\le \log H$.
%%%%%%%%%%%%%%%%%%
\bigskip

We may deduce distributional inequalities from this  hypothesis.
%%%%%%%%%%%%%%
\begin{lem}
  \label{Distrib}
  Under the hypothesis $(Hyp.)$, we have
  \begin{equation*}
    \biggl|\Bigl\{
    n\in \mathcal{N}:\ \Bigl|\sum_{x\in\mathcal{X}}c(x)e(nx)\Bigr|
    \ge \lambda \sqrt{\sum_{x\in\mathcal{X}}|c(x)|^2}
    \Bigr\}\biggr|
    \le 4 H
    e^{-\lambda^2/4}.
  \end{equation*}
\end{lem}
%%%%%%%%%%%%%%
\noindent
See the second half of \cite[Lemma 4.33]{Tao-Vu*06} by T.~Tao \& V.~Vu.
As these inequalities may be compared with the classical Chernoff bound
for sums of independent random variables, one may say that this shows
the $e(nx)$ behaves ``independently''.

%%%%%%%%
\begin{proof}
  Let us use Hypothesis~$(Hyp.)$ on $\tilde{c}(x)=e(\theta)\sigma\,c(x)$
  for a positive parameter $\sigma$ and a phase $\theta$ that we may
  choose. Set $C=\sum_{x}|c(x)|^2$. We find that
  \begin{align*}
    \Bigl|\bigl\{
    n\in\mathcal{N}, \sigma \Re e(\theta)\sum_{x}c(x)e(nx)\ge \sigma \lambda \sqrt{C}
    \bigr\}\Bigr|\,\exp(\sigma \lambda \sqrt{C})
    \le H
    \exp(\sigma^2 C/2).
  \end{align*}
  On selecting $\sigma=\lambda/\sqrt{C}$, we reach
  \begin{equation}
    \label{Ibase}
    \bigl|\bigl\{
    n\le \mathcal{N}, \Re e(\theta)\sum_{x}c(x)e(nx)\ge \lambda \sqrt{C}
    \bigr\}\bigr|
    \le H
    e^{-\lambda^2 /2}.
  \end{equation}
  To infer inequalities on the modulus, let us notice that
  \begin{equation*}
    \max \Bigl(\bigl|\Re \sum_{x}c(x)e(nx)\bigr|
    ,\bigl|\Im \sum_{x}c(x)e(nx)\bigr|\Bigr)\sqrt{2}\ge
     \bigl|\sum_{x}c(x)e(nx)\bigr|.
   \end{equation*}
   Therefore
   \begin{align*}
     (1/\sqrt{2})
    \Bigl|\sum_{x}c(x)e(nx)\Bigr|
     \le
     \max \Bigl(
     &\Re \sum_{x}c(x)e(nx),
     \Re (-1)\sum_{x}c(x)e(nx),
     \\&\Re i\sum_{x}c(x)e(nx),
     \Re (-i)\sum_{x}c(x)e(nx)
     \Bigr),
   \end{align*}
   so inequality~\eqref{Ibase} applies (which $\lambda/\sqrt{2}$
   rather than $\lambda$). The lemma is proved.
\end{proof}
%%%%%%%%

%%%%%%%%%%%%%
\begin{lem}
  \label{AS}
  When $x>0$, there exists $\theta\in(0,1)$ such that
  \begin{equation*}
    |\Gamma(x+1)|= \sqrt{2\pi}
    |x|^{x+\frac12}\exp\biggl(-x+\frac{\theta}{12x}\biggr).
  \end{equation*}
\end{lem}
%%%%%%%%%%%%%
This is from the book \cite[(6.1.38)]{Abramowitz-Stegun*64} by {M.}
Abramowitz and {I.A.} Stegun.

We now turn towards L${}^p$-inequalities. The next one is the
very similar, when replacing the summation over a finite abelian group
by $n\le N$, of \cite[Theorem 1.1]{Shkredov*11} by I.~Shkredov, which this author calls
Rudin's inequality.
%%%%%%%%%%%%%%
\begin{lem}
  \label{LpQRI}
  Under the hypothesis $(Hyp.)$ and for any real number $p\ge1$, we have
  \begin{equation*}
    \sum_{n\in\mathcal{N}}
    \Bigl|\sum_{x\in\mathcal{X}}c(x)e(nx)\Bigr|^p
    \le
    4 (\tfrac95\sqrt{p})^{p}H
    \,\biggl(\sum_{x\in\mathcal{X}}|c(x)|^2\biggr)^{p/2}.
  \end{equation*}
\end{lem}
%%%%%%%%%%%%%%
In \cite[Proposition 1]{Green*04}, B.~Green proves a resembling inequality on
$\Z{N}$, with the slightly worse constant $(12\sqrt{p})^{p}$.
%%%%%%%%%% 
\begin{proof}
  Set $f(n)=\sum_{x}c(x)e(nx)$ and $F=\sqrt{\sum_{x}|c(x)|^2}$. We
  define
  \begin{equation*}
    \mathcal{N}(\lambda)
    =\bigl\{n\in\mathcal{N}:|f(n)|\ge\lambda F\bigr\}
  \end{equation*}
  as well
  as $N(\lambda)=|\mathcal{N}(\lambda)|$. Let $0\le
  \lambda_1<\lambda_2<\cdots<\lambda_K$ be the sequence of values
  taken by $|f(n)|/F$, and we set $N(\lambda_{K+1})=0$ as well as $\lambda_0=0$.
  We readily find that
  \begin{align*}
    \sum_{n\in\mathcal{ N}}|f(n)|^p/F^p
    &=\sum_{k\le K} \lambda_k^p
    \bigl(N(\lambda_k)-N(\lambda_{k+1})\bigr)
    =\sum_{k\le K} N(\lambda_k)(\lambda_k^p-\lambda_{k-1}^p)
    \\&=\sum_{k\le K} N(\lambda_k)\int_{\lambda_{k-1}}^{\lambda_k}pt^{p-1}dt
    = \int_{0}^{\infty}pN(t)t^{p-1}dt.
  \end{align*}
  We may apply Lemma~\ref{Distrib} to majorize $N(t)$.
  We are thus led to find an upper bound for
  $\int_{0}^{\infty}pe^{-t^2/4}t^{p-1}dt$. Set $u=t^2/4$, so that
  $dt=du/\sqrt{u}$, and get that
  \begin{align*}
    \int_{0}^{\infty}pe^{-t^2/4}t^{p-1}dt
    &=
    2^{p-1}p\int_{0}^{\infty}u^{p/2}e^{-u}du/u
    =p2^{p-1}\Gamma(p/2)
    \\&=
    2^{p}\Gamma\Bigl(1+\frac{p}{2}\Bigr)
    \le \sqrt{2\pi}2^{p}(p/2)^{\frac{p+1}{2}}e^{-p/2}e^{1/(6p)}
    \\&\le \sqrt{\pi}2^{p/2} p^{\frac{p+1}{2}}e^{-p/2}e^{1/(6p)}
    \le (\tfrac95\sqrt{p})^p
  \end{align*}
  by Pari/GP.
  % g(p)=log(sqrt(Pi)*exp((log(2)-1)*p/2+1/6/p)*sqrt(p))/p
  % gg(p)=2/sqrt(p)*gamma(1+p/2)^(1/p)
\end{proof}
%%%%%%%%%% 
%%%%%%%%%%%%%%
\begin{thm}
  \label{RA}
  Under the hypothesis $(Hyp.)$, let $f$ be a function over
  $\mathcal{N}$ with support
  inside a given subset~$S$. We have
  \begin{equation*}
    \sum_{x\in\Xcal}
    \biggl|\sum_{n\in\mathcal{N}}f(n)e(nx)\biggr|^2
    \le
    9
    \, |S|\, \|f\|_2^2\, \log\frac{8H}{|S|}
    .
  \end{equation*}
  More generally, for any real number $\ell\ge1$, we have
  \begin{equation*}
    \biggl(\sum_{x\in\Xcal}
    \biggl|\sum_{n\in\mathcal{N}}f(n)e(nx)\biggr|^{\ell+1}\biggr)^2
    \le
    9\,
    |S|\,\|f\|_2^2\,\sum_{x\in\Xcal}
    \biggl|\sum_{n\in\mathcal{N}}f(n)e(nx)\biggr|^{2\ell}
    \log\frac{8H}{|S|}
    .
  \end{equation*}
\end{thm}
%%%%%%%%%%%%%%

%%%%%%%%%% 
\begin{proof}
  Let us start, with an obvious notation, from
  \begin{align*}
    \biggl|\sum_x\sum_n f(n)e(nx) g(x)\biggr|
    &\le
      \Bigl(\sum_n |f(n)|^q\Bigr)^{1/q}
      \biggl(\sum_n \Bigl|\sum_{x}g(x)e(nx)\Bigr|^p\biggr)^{1/p}
    % \\&\le
    %   \Bigl(\sum_n |f(n)|^q\Bigr)^{1/q}
    %   \biggl(\sum_n \Bigl|\sum_{x}g(x)e(nx)\Bigr|^p\biggr)^{1/p}
    \\&\le
      \Bigl(\sum_n |f(n)|^q\Bigr)^{1/q}
    \biggl(4(\tfrac95\sqrt{p})^pH
    \Bigl(\sum_x |g(x)|^2\Bigr)^{p/2}
    \biggr)^{1/p}
  \end{align*}
  where Lemma~\ref{LpQRI} was used in the second step.
  We select $g(x)=\sum_n \overline{f(n)}e(-nx)$ to obtain
  \begin{equation*}
    \sum_{x\in\Xcal}
    \Bigl|\sum_{n\in\mathcal{N}}f(n)e(nx)\Bigr|^2
    \le \bigl(\tfrac95\bigr)^2 \Bigl(\sum_n |f(n)|^q\Bigr)^{2/q}
    p\, H^{2/p}.
  \end{equation*}
  This is the dual form of Lemma~\ref{LpQRI}.  Let us now use the
  hypothesis on the support of $f$. We assume that $p> 2$, i.e. $q< 2$.
  Let us again us the H\"older inequality, this time with exponents $a=2/q$ and $b$
  defined $1/b=1-q/2$. We find that
  \begin{equation*}
    \Bigl(\sum_n |f(n)|^q\Bigr)^{2/q}
    \le\sum_n |f(n)|^2 |S|^{1-\frac2p}.
  \end{equation*}
  Therefore
  \begin{equation*}
    \sum_{x\in\Xcal}
    \Bigl|\sum_{n\in\mathcal{N}}f(n)e(nx)\Bigr|^2
    \le 
    \bigl(\tfrac95\bigr)^2 p
    (H/|S|)^{2/p}
    |S|\sum_n |f(n)|^2.
  \end{equation*}
  The choice $p = \log(8 H/|S|)\ge 2$ gives the first inequality. We
  may also take $g(x)=\sum_n \overline{f(n)}e(-nx)|\sum_n
  f(n)e(nx)|^{\ell-1}$ for some $\ell\ge1$. Then we use
  \begin{equation*}
    \Bigl(\sum_x |g(x)|^2\Bigr)^{1/2}
    =
    \biggl(\sum_x \biggl|\sum_n
    f(n)e(nx)\biggr|^{2\ell}\biggr)^{1/2}.
  \end{equation*}
  This ends the proof of our theorem.
\end{proof}
%%%%%%%%%%

%%%%%%%%%%%%%%%%%%
\section{A quantitative version of Rudin's inequality}
%%%%%%%%%%%%%%%%%%
\label{Integers}
The handling of the integers from $[1,N]$ will almost exclusively on the
next lemma. It is due to A. Selberg, see \cite[Section
20]{Selberg*91} or
\cite{Vaaler*85} by J.D. Vaaler.
%%%%%%%%%%
\begin{lem}
  \label{BeurlingF}
  Let $M\in\mathbb{R}$, and $N$ and $\delta$ be positive real number.
  There exists a smooth function $\psi$ on $\mathbb{R}$ such that
  \begin{itemize}
  \item The function $\psi$ is non-negative.
  \item When $t\in [M,M+N]$, we have $\psi(t)\ge 1$.
  \item $\psi(0)=N+\delta^{-1}$.
  \item When $|\alpha|>\delta$, we have $\hat{\psi}(\alpha)=0$.
  \item We have $\psi(t)=\Ocal_{M,N,\delta}(1/(1+|t|^2))$.
  \end{itemize}
\end{lem}
%%%%%%%%%%

%%%%%%%%%%%%
\begin{lem}
  \label{QRI}
  Under the hypotheses of Theorem~\ref{LpQRI}, we have
  \begin{equation*}
    \log \sum_{n\le N}\biggl|\exp\biggl(
    \sum_{x\in \Xcal}c(x) e(xn)\biggr)\biggr|
    \le \tfrac12 \sum_{x\in\Xcal}|c(x)|^2+\log(N+\delta_\star^{-1}).
  \end{equation*}
\end{lem}
%%%%%%%%%%%%
Let us specify that the summation is over positive~$n$.
We follow the proof of Lemma 4.33 in \cite{Tao-Vu*06} by T. Tao \&
V.H. Vu.
See also \cite[Proposition 17]{Green*02a} by B. Green.
There are $2^{|\Xcal|}$ subsums $\sum_{x\in\Acal}x$. As they are
suposed to be all distincts in $\mathbb{R}/\mathbb{Z}$,
the minimal spacement is at most $2^{-|\Xcal|-1}$. 
%%%%%%%%
\begin{proof}
  Let us consider
  \begin{equation}
    \label{defSigma}
    \Sigma=\sum_{n\le N}\biggl|\exp\biggl(
    \sum_{x\in \Xcal}c(x) e(xn)\biggr)\biggr|
    =
    \sum_{n\le N}\exp\biggl(\Re\,
    \sum_{x\in \Xcal}c(x) e(xn)\biggr).
  \end{equation}
  Let us write $c(x)=|c(x)|e(\theta_x)$. As
  \begin{equation}
    \label{eq:2}
    \forall y\ge0, \forall t\in[-1,1],\quad e^{ty}\le \ch y+t\sh y,
  \end{equation}
  we find that
  \begin{equation}
    \label{eq:1}
    \exp\bigl(|c(x)|\Re e(xn+\theta_x)\bigr)
    \le \ch(|c(x)|) + \sh(|c(x)|)\Re e(xn+\theta_x).
  \end{equation}
  We gather these inequalities to reach that $\Sigma$ is bounded
  above by
  \begin{equation}
    \label{step1}
     \sum_{n\in\mathbb{Z}}\psi(n)
    \prod_{x\in \Xcal}\Bigl(
    \ch(|c(x)|)
    + \tfrac12\sh(|c(x)|) e(xn+\theta_x)
    + \tfrac12\sh(|c(x)|) e(-xn-\theta_x)
    \Bigr).
  \end{equation}
  Let $(\Acal,\Bcal,\Ccal)$ range over partitions of $\Xcal$. By using
  (compare Taylor series)
  \begin{equation*}
    \ch y\le e^{y^2/2},
  \end{equation*}
  we reach
  \begin{align*}
    \Sigma
    &\le
    \Bigl(N+
    \sum_{(\Acal,\Bcal,\Ccal)\neq(\emptyset,\emptyset,\Xcal)}\delta_\star^{-1}2^{-|\Acal|-|\Bcal|}
    \Bigr)e^{\sum_{x\in\Xcal}|c(x)|^2/2}
    \\&\le
    \bigl(N+
    \delta_\star^{-1}2^{|\Xcal|}
    \bigr)e^{\sum_{x\in\Xcal}|c(x)|^2/2}.
  \end{align*}
  Indeed, we had to sum $\sum_{n\in\mathbb{Z}}\psi(n)e(ny)$ for
  $y=\sum_{x\in\Acal}x-\sum_{x\in\Bcal}x\neq0$. 
\end{proof}
%%%%%%%%

%%%%%%%%%% 
\begin{proof}[Proof of Theorem~\ref{LpQRI*}]
  Lemma~\ref{QRI} is prepared to apply Theorem~\ref{RA} which is
  precisely what we claimed.
\end{proof}
%%%%%%%%%%

%%%%%%%%%%%%%%%%%%
\section{Rudin's inequality on primes}
%%%%%%%%%%%%%%%%%%
\label{Primes}
We shall handle the primality condition via an enveloping sieve, as
initially in~\cite{Ramare*95}, but with an additional parameter
introduced in~\cite[Section 4, (8)]{Ramare*22}. This is more fully
described in the future paper~\cite[Section 3]{Ramare*25-1}.
As a matter of fact, the reader may proceed without this paramater. It
is included only for some applications we have in mind where it will be
required. The proof is not more difficult.

We proceed as in the case of integers and establish the following inequality.
%%%%%%%%%%%%
\begin{lem}
  \label{QRIp}
  Under the hypothesis and notation of Theorem~\ref{LpQRI*p}, we have
  \begin{equation*}
    \sum_{N^\kappa<p\le N}\biggl|\exp\biggl(
    \sum_{x\in \Xcal}c(x) e(xp)\biggr)\biggr|
    \le 8\frac{N+\delta_*^{-1}(\sqrt{N},z_0)}{\kappa\log N}(\log z_0)
    e^{\frac12 \sum_{x\in\Xcal}|c(x)|^2}.
  \end{equation*}
\end{lem}
%%%%%%%%%%%%

%%%%%%%%%%% 
\begin{proof}
  We proceed as in the proof of Lemma~\ref{QRI} until we reach~\eqref{step1}.
  There we bound above the characteristic
  function of the primes by the enveloping sieve of \cite[Section
  4]{Ramare*22} with $z=N^{\kappa/2}$.
  This leads to
  \begin{multline}
    \sum_{(\Acal,\Bcal,\Ccal)}
    \prod_{x_a\in\Acal}\frac{|c(x_a)|}{2}
    \prod_{x_b\in\Bcal}\frac{|c(x_b)|}{2}
    \prod_{x_c\in\Ccal}|c(x_c)|
    \\
    \sum_{n\in\mathbb{Z}}\beta_{z_0,z}(n)
    e\bigl((x_{\mathcal{A}}-x_{\mathcal{B}})n+\theta_{\mathcal{A}}-\theta_{\mathcal{B}}\bigr)
    \psi(n)
    \label{begsieve}
  \end{multline}
  where $x_{\mathcal{A}}=\sum_{a\in\mathcal{A}}x_a$, similarly for
  $x_{\mathcal{B}}$ and
  $\theta_{\mathcal{A}}=\sum_{a\in\mathcal{A}}\theta_a$, and similarly
  for $\theta_{\mathcal{B}}$, with the notation of~\eqref{eq:1} and~\eqref{step1}.
  Since $\beta_{z_0,z}(n)=(\sum_{d|n}\lambda_d)^2$, this may also be rewritten as
  \begin{multline*}
    \sum_{(\Acal,\Bcal,\Ccal)}
    \prod_{x_a\in\Acal}\frac{|c(x_a)|}{2}
    \prod_{x_b\in\Bcal}\frac{|c(x_b)|}{2}
    \prod_{x_c\in\Ccal}|c(x_c)|
    \\
    \sum_{d_1,d_2\le z}\lambda_{d_1}\lambda_{d_2}
    \sum_{m\in\mathbb{Z}}
    e\bigl((x_{\mathcal{A}}-x_{\mathcal{B}})m[d_1,d_2]+\theta_{\mathcal{A}}-\theta_{\mathcal{B}}\bigr)
    \psi([d_1,d_2]m).
  \end{multline*}
  Poisson summation formula enables us to rewrite the sum over $m$ in
  the form
  \begin{equation*}
    \frac{e(\theta_{\mathcal{A}}-\theta_{\mathcal{B}})}{[d_1,d_2]}
    \sum_{k\in\mathbb{Z}}
    \hat{\psi}\biggl(\frac{k-(x_{\mathcal{A}}-x_{\mathcal{B}})[d_1,d_2]}{[d_1,d_2]}\biggr).
  \end{equation*}
  This shows that this sum sum vanishes when
  $x_{\mathcal{A}}-x_{\mathcal{B}}\neq 0$ and equals
  \begin{equation}
    \label{maintsievestep}
    (N+\delta_*^{-1}(z_0,z))
    \sum_{d_1,d_2\le z}\frac{\lambda_{d_1}\lambda_{d_2}}{[d_1,d_2]}
    =\frac{N+\delta_*^{-1}(z,z_0)}{G(z;z_0)}
  \end{equation}
  otherwise. We conclude the proof as we did for the one of
  Lemma~\ref{QRI}. \cite[Lemma 2.6]{Ramare*22} provides us with a
  lower bound for $G(z;z_0)$. We thus reach the upper bound
  \begin{equation*}
    \biggl(e^\gamma \frac{N+\delta_*^{-1}(z,z_0)}{\log z}\log (2z_0)
    \biggr)e^{\frac12\sum_{x}|c(x)|^2}.
  \end{equation*}
  We  notice that $2e^\gamma\frac{N}{\log N}\log
  (2z_0)\le 8 \frac{N}{\log N}\log z_0$. The proof is complete.
\end{proof}
%%%%%%%%%%%

%%%%%%%%
\begin{proof}[Proof of Theorem~\ref{LpQRI*p}]
  Lemma~\ref{QRIp} is prepared to apply Theorem~\ref{RA}. We obtain
  the inequality
  \begin{equation*}
    \sum_{x\in\Xcal}
    \biggl|\sum_{p\le N}f(p)e(px)\biggr|^2
    \le
    9
    \, |S|\, \|f\|_2^2\,
    \log\biggl(32\frac{N+\delta_*^{-1}(\sqrt{N},z_0)}{\kappa|S|\log N}\log z_0
    \biggr).
  \end{equation*}
  Theorem~\ref{LpQRI*p} follows swiftly from there.
\end{proof}
%%%%%%%%

%%%%%%%%%
\begin{proof}[Proof of Corollary~\ref{cuspsmodU}]
  We use Theorem~\ref{LpQRI*p} with $z_0=2$ and $\kappa=1/4$. We have
  $\delta_*(N^\kappa,z_0)>1/(N^{1/4}U_1)\ge 1/U$. Therefore, with $K=\frac{N}{|S|\log
    N}$, we find that
  \begin{equation*}
    c\le 9\exp\frac{\log\Bigl(64\frac{2UK}{N}\log 2\Bigr)}{\log 8K}
    \le 30e^{U/N}.
  \end{equation*}
  %9*exp(log(64*2*log(2)/8)/log(8))
  The proof is complete.
\end{proof}
%%%%%%%%%

%%%%%%%%%%%%%%%%%%%%%%%%%%%%%%%
\section{Proof of Theorem~\ref{StrutCusps}} 
%%%%%%%%%%%%%%%%%%%%%%%%%%%%%%%
%%%%%%%%%%%%
\begin{proof}[Proof of Theorem~\ref{StrutCusps}]
  Indeed, for $i\in\{1,2\}$, take for $\mathcal{D}_i=\{a_1(i),\cdots,a_{K(i)}(i)\}$ a maximal
  set inside $\mathscr{C}/U_i\mathbb{Z}\subset \Z{U_i}$ that is \emph{dissociate},
  i.e. such that all its sumset sums
  are distinct. Consider the lift of
  $\mathcal{E}=\mathcal{D_1}\times\mathcal{D}_2$ by the Chinese
  Remainder Theorem map $\varphi:\Z{U_1U_2}\rightarrow\Z{U_1}\times\Z{U_2}$.
  For any $a_k(2)$, the set
  $\varphi^{-1}(\mathcal{D_1}\times\{a_k(2)\})$ satisfies the
  assumptions of Corollary~\ref{cuspsmodU}, so its cardinality $K(1)$
  is bounded by $30A^2\log(8K)e^{U_1U_2/N}$.
  %Let us set $U_2=U/U_1$
  Any further point $x'\in\mathcal{C}$ has that, for $i\in\{1,2\}$,
  the set 
  $\{x',a_1(i),\cdots,a_{K(i)}(i)\}/U_i\mathbb{Z}$ is not dissociate, from which the reader
  will readily conclude that there exists $\Acal,\Bcal\in\mathcal{D}_i$ such
  that
  \begin{equation*}
    x'+\sum_{x\in\Acal}x=\sum_{x\in\Bcal}x.
  \end{equation*}
  The lemma follows swiftly from there. 
\end{proof}
%%%%%%%%%%%%

%%%%%%%%%%%%%%%%%%
\section{On the Selberg sieve for squares}
%%%%%%%%%%%%%%%%%%
\label{SSQ}
How to use the Selberg sieve for squares is completely described in
the conventional setting in~\cite[Chapter 11]{Ramare*06}.
The involved quantities are used and evaluated already
in~\cite[Theorem 5.4]{Ramare*06}.
%%%%%%%%%%%%%%%%%%%%%
\subsection*{Pointwise upper bound}
%%%%%%%%%%%%%%%%%%%%%
We define, for any positive integer $q$, the set
$\mathcal{K}_q\subset\Z{q}$ ro be the set of squares modulo~$q$. For
simplicity we restrict our attention to moduli~$q$ that are odd and
square-free. The function $q\mapsto|\mathcal{K}_q|$ is multiplicative
and takes the value $(p+1)/2$ at the prime~$p$. We seek real
parameters $(\lambda^\sharp_q)_{q\le Q}$ that are such that
$\sum_q\lambda^\sharp_q=1$ and whose support in on odd square-free
integers. We chose them so as to optimize the main term in the upper
bound:
\begin{equation}
  \label{eq:4}
  \sum_{\substack{n\le N\\ \text{$n$ is a square}}}1
  \le
  \sum_{\substack{n\le N}}\biggl(
  \sum_{\substack{q\le z\\ n\in\mathcal{K}_q}}\lambda^\sharp_q\biggr)^2.
\end{equation}
This main term is
\begin{equation}
  \label{eq:5}
  \text{MT}
  =N\sum_{q_1,q_2\le Q}\frac{|\mathcal{K}_{[q_1,q_2]}|}{[q_1,q_2]}
  \lambda^\sharp_{q_1}\lambda^\sharp_{q_2}
  =
  N\sum_{q_1,q_2\le Q}\frac{(q_1,q_2)}{|\mathcal{K}_{(q_1,q_2)}|}
  \frac{|\mathcal{K}_{q_1}|\lambda^\sharp_{q_1}}{q_1}
  \frac{|\mathcal{K}_{q_2}|\lambda^\sharp_{q_2}}{q_2}.
\end{equation}
We introduce the multiplicative function $h$ which vanishes on
non square-free integers and whose value on primes is given by
$h(2)=0$ and $h(p)=(p-1)/(p+1)$ when $p$ is an odd prime. This
function satisfies the convolution identity
\begin{equation*}
  \frac{q}{|\mathcal{K}_{q}|}
  =\sum_{d|q}h(q)
\end{equation*}
when $q$ is odd and square-free. Therefore
\begin{equation*}
  \text{MT}
  =
  N\sum_{d\le Q}h(d)\biggl(\sum_{\substack{q\le z\\ d|q}}
  \frac{|\mathcal{K}_{q}|\lambda^\sharp_{q}}{q}\biggr)^2.
\end{equation*}
We are thus led, as in the classical Selberg approach, to a quadratic
optimization problem under a linear constraint. We cut the story short
and select
\begin{equation}
  \label{deflambdasharpq}
  \lambda^\sharp_{q}=\frac{q}{|\mathcal{K}_{q}|}
  \sum_{\substack{d\le z\\ q|d}}\mu(d/q)/G^\sharp(z),
  \quad
  G^\sharp(z)=\sum_{q\le z}h(\delta).
\end{equation}
We simply check that
\begin{align*}
  G^\sharp(z)\sum_{q\le z}\lambda^\sharp_{q}
  &=\sum_{q\le z}(\1\star h)(q)
  \sum_{\substack{d\le z\\ q|d}}\mu(d/q)
  \\&=
  \sum_{\substack{d\le z}}(\mu\star\1\star h)(q)=G^\sharp(Q).
\end{align*}
The coefficient $(\lambda^\sharp_q)$ are thus very similar to the
Selberg coefficients $(\lambda_d)$ except that are \emph{not}
bounded. How to recover the Selberg coefficients from these is
described fully in~\cite[Chapter 11]{Ramare*06} and in a paragraph
below. Both may be skipped here as we shall only use the density
$G^\sharp(z)$.

%%%%%%%%%%%%%%%%%%%%%
\subsection*{Fourier decomposition of the sieve}
%%%%%%%%%%%%%%%%%%%%%
The sets $\mathcal{K}_d$ are subsets of $\Z{d}$, which is not apparent
in the writing. We now propose a decomposition that makes this point
foremost. Here is how it goes:
\begin{equation*}
  \beta_z(n)
  =\biggl(\sum_{q\le z}\lambda_q^\sharp\1_{\mathcal{K}_q}(n)\biggr)^2
  =\sum_{q_1,q_2\le z}\lambda_{q_1}^\sharp\lambda_{q_2}^\sharp\1_{\mathcal{K}_{q_1}}(n)
  \1_{\mathcal{K}_{q_2}}(n).
\end{equation*}
The next step is to notice that
$\1_{\mathcal{K}_{q_1}}\1_{\mathcal{K}_{q_2}}=\1_{\mathcal{K}_{[q_1,q_2]}}$
where $[q_1,q_2]$ is the lcm of~$q_1$ and $q_2$. We develop this
function in Fourier series
\begin{equation*}
  \1_{\mathcal{K}_{q}}(n)=\frac{1}{q}\sum_{a\mod q}\eta(q;a)e(na/q)
\end{equation*}
for some coefficients $\eta(q;a)$. We next introduce \emph{primitive}
additive characters via
\begin{equation*}
  \1_{\mathcal{K}_{q}}(n)=\sum_{d|q}\sum_{a\mode d}\eta(q;bq/d)e(nb/d)
\end{equation*}
where $b\mode d$ denotes a sum over integers $b$ in $\{1,\cdots,d\}$
that are coprime to~$d$. After some shuffling, we reach
\begin{equation}
  \label{Fourierbeta}
  \beta_z(n)
  =\sum_{d\le z^2}\sum_{b\mode d}w_d(\mathcal{K}, b/d)e(na/d),
\end{equation}
with
\begin{equation}
  \label{defw}
  w_d(\mathcal{K}, b/d)
  =
  \sum_{\substack{q_1,q_2\le z\\ d|[q_1,q_2]}}
  \frac{\lambda_{q_1}^\sharp\lambda_{q_2}^\sharp}{[q_1q_2]}
  \eta([q_1,q_2];b[q_1,q_2]/d).
\end{equation}
One could study there coefficients much more but their sole existence
will be enough below. We call Eq.~\eqref{Fourierbeta} the
\emph{Fourier expansion} of $\beta$, as in \cite[End of subsection
4.1]{Ramare-Ruzsa*01}. The weights $\beta_{z_0,z}(n)$ would be defined
in just that same manner but carrying a condition $(q,P(z_0))=1$
throughout. This condition is trivial when $z_0=2$.

%%%%%%%%%%%%%%%%%%%%%
\subsection*{Recovering the usual Selberg coefficients}
%%%%%%%%%%%%%%%%%%%%%
The material of this subsection is not used in this paper, it is
only included for the readers to see properly the connections between
different viewpoints.

A classical presentation of the sieve, as for instance in the book
\cite{Bombieri*74} by E.~Bombieri, starts from sets
$\Omega_p\subset\Z{p}$ that are to be \emph{avoided}. We then look
typicall at the integers $n\le N$ that are so that $n\notin\Omega_p$
for every $p$ into consideration. Reverting to the above notation, we
have $\Omega_p=\Z{p}\setminus\mathcal{K}_p$. In~\cite{Ramare*06}, we
extended the notion of $(\Omega_p)$ to $(\Omega_q)$, where $q$ is any
modulus, called this one the \emph{bordering sytem} and in fact denoted
it by $(\mathcal{L}_q)$. But let us stick in this short presentation
to the notation $\Omega_q$. The set
$\Omega_q$ is simply obtained\footnote{This is because we restrict our
  attention to square-free moduli~$q$'s.} by glueing together the $\Omega_p$ for
$p|q$ via the Chinese remainder Theorem, but $n\notin\Omega_p$ for $p|q$
\emph{does not translate} in $n\notin\Omega_q$. We may solve
this difficulty by resorting to indicator functions. We have
$\1_{\mathcal{K}_p}=\1-\1_{\Omega_p}$ and therefore, when $d$ is square-free:
\begin{equation}
  \label{eq:8}
  \1_{\mathcal{K}_d}=\sum_{\ell|d}\mu(\ell)\1_{\Omega_\ell}.
\end{equation}
We continue and infer that
\begin{equation*}
  \sum_{\substack{q:n\in\mathcal{K}_q}}\lambda_q^\sharp
  =
  \sum_{\substack{q}}\lambda_q^\sharp\1_{\mathcal{K}_q}(n)
  =
  \sum_{\substack{\ell}}\lambda_\ell\1_{\Omega_\ell}(n)
\end{equation*}
where $\lambda_\ell=\mu(\ell)\sum_{q:\ell|q}\lambda_q^\sharp$, so that
\begin{equation}
  \label{eq:10}
  G^\sharp(z)\lambda _\ell
  =\mu(\ell)\sum_{q:\ell|q}(\1\star h)(q)
  \sum_{\substack{d\le z\\ q|d}}\mu(d/q)
  =\mu(\ell)(1\star h)(\ell)
  \sum_{\substack{m\le z/\ell\\ (m,\ell)=1}}h(m).
\end{equation}
These coefficients $(\lambda_\ell)$ are exactly the ones that we
obtain in the classical presentation of the Selberg sieve, the
condition $\sum_{q}\lambda^\sharp_q=1$ translating in $\lambda_1=1$. Yet again
the method of \cite{van-Lint-Richert*65} by J.E.~van Lint and
H.E.~Richert may be adapted to show that $|\lambda_\ell|\le 1$, but
the situation changes dramatically here because we in fact have (on
anticipating on the evaluation of $G^\sharp(z)$ in the next paragraph)
\begin{equation}
  \label{eq:11}
  \lambda_\ell\ll 2^{\omega(\ell)}/\ell,
\end{equation}
as simple consequence of the bound $\sum_{\substack{m\le
    z/\ell\\ (m,\ell)=1}}h(m)\le z/\ell$.

%%%%%%%%%%%%%%%%%%%%% 
\subsection*{Density evaluation}
%%%%%%%%%%%%%%%%%%%%%
We want to sieve up $z=\sqrt{N}$ and the main quantity
to be evaluated (minorized is enough) is:
\begin{equation}
  \label{eq:3}
  G^\sharp(z)=\sum_{q\le z}h(\delta)
  =\sum_{\substack{q\le z\\ (q,2)=1}}\mu^2(q)\prod_{p|q}\frac{p-1}{p+1}
\end{equation}
denoted by $G_{d}(z)$ at the bottom of \cite{Ramare*06}, for
$d=1$. The function $h$ therein is indeed the one we have defined
above. The estimation given page~49 with $f=1$ and $u=0$, so the function $a$ reduces to
the function~$h$, gives us
\begin{equation}
  \label{minGshap}
  G^\sharp(z) = B(1)\bigl(z +\Ocal^*(6.7\, z^{3/4})\bigr)
  \ge 0.35 (1-6.7/z^{1/4}) z\ge 0.326\,z
\end{equation}
when $z\ge 10^8$. A Pari/GP script rapidly shows that $G^\sharp(z)/z\ge 0.304$
when $z\le 10^9$, the minimum being reached just before~$z=179$. We
will use the bound $G^\sharp(z)/z\ge 0.304\ge 1/8$ for $z\ge1$ so as
to get numerics exactly similar to the ones in the case of the primes.
%%%%%%%%%%%%%%%%%%
\section{Rudin's inequality for squares}
%%%%%%%%%%%%%%%%%%
The Selberg sieve for squares is briefly described in Section~\ref{SSQ}.
%%%%%%%%%%%%%%%%%%
\begin{lem}
  \label{QRIQ}
  Under the hypothesis and notation of Theorem~\ref{LpQRI*p}, we have that
  \begin{equation*}
    \sum_{n\le \sqrt{N}}\biggl|\exp\biggl(
    \sum_{x\in \Xcal}c(x) e(xn^2)\biggr)\biggr|
    \le 8\frac{N+\delta_*^{-1}(\sqrt{N},2)}{\sqrt{N}}
    e^{\frac12 \sum_{x\in\Xcal}|c(x)|^2}.
  \end{equation*}
\end{lem}
%%%%%%%%%%%%%%%%%%

%%%%%%%%%%%%%
\begin{proof}
  The proof follows faithfully the one of Theorem~\ref{QRIp} until
  Eq.~\eqref{begsieve}, where we dispense of the parameter $z_0$ and
  employ the decomposition given in~\eqref{Fourierbeta}.
  Our upper bound becomes
   \begin{multline*}
    \sum_{(\Acal,\Bcal,\Ccal)}
    \prod_{x_a\in\Acal}\frac{|c(x_a)|}{2}
    \prod_{x_b\in\Bcal}\frac{|c(x_b)|}{2}
    \prod_{x_c\in\Ccal}|c(x_c)|
    \\
    \sum_{d\le z^2}\sum_{b\mode d}w_d(\mathcal{K}, b/d)
    \sum_{n\in\mathbb{Z}}
    e\biggl((x_{\mathcal{A}}-x_{\mathcal{B}})n+\theta_{\mathcal{A}}-\theta_{\mathcal{B}}
    +\frac{bn}{d}\biggr)
    \psi(n).
  \end{multline*}
  Yet again, Poisson summation formula enables us to rewrite the sum over $n$ in
  the form
  \begin{equation*}
    e(\theta_{\mathcal{A}}-\theta_{\mathcal{B}})
    \sum_{k\in\mathbb{Z}}
    \hat{\psi}\biggl(\frac{dk-b-(x_{\mathcal{A}}-x_{\mathcal{B}})d}{d}\biggr).
  \end{equation*}
  At this level, we simply need to put
  another evaluation for the sieve density $G(z;z_0)$, and replace it
  by~\eqref{minGshap} for $Q=\sqrt{N}$. The lemma follows readily.
\end{proof}
%%%%%%%%%%%%%

%%%%%%%%%%
\begin{proof}[Proof of Theorem~\ref{LpQRI*Q}]
  Lemma~\ref{QRIQ} is prepared to apply Theorem~\ref{RA}. As it is
  exactly similar to Lemma~\ref{QRIp}, the reader
  will complete the proof without any difficulty.
\end{proof}
%%%%%%%%%%

%\bibliographystyle{authordate1}
%\bibliography{Local.bib}
%\bibliographystyle{plain}
%\printbibliography

\end{document}